\documentclass[11pt]{article}
\usepackage{amssymb,amsfonts,amsmath,latexsym, epsfig,mathrsfs,colordvi}
\usepackage{graphicx}
\newtheorem{theo}{Theorem}

\newtheorem{exam}[theo]{Example}
\newtheorem{lem} [theo]{Lemma}

\newcommand{\Star}{\mathrm{Star}}

\makeatletter \@addtoreset{equation}{section}
\@addtoreset{theo}{}\makeatother

\def\qed{\hfill \rule{4pt}{7pt}}
\def\pf{\noindent {\it Proof.} }

\begin{document}

\title{Counting Labelled Trees with Given Indegree
Sequence}

\author{Rosena R. X. Du\footnote{Department of Mathematics, East China Normal University,
Shanghai 200241, P. R. China, Email: rxdu@math.ecnu.edu.cn. }\ \ and
Jingbin Yin\footnote{Department of Mathematics, Massachusetts
Institute of Technology, Cambridge, MA 02139, USA, Email:
jbyin@math.mit.edu.}}

\date{Mar 27, 2008}
\maketitle

\vskip 0.7cm \noindent{\bf Abstract.} For a labelled tree on the
vertex set $[n]:=\{1,2,\ldots, n\}$, define the direction of each
edge $ij$ to be $i\rightarrow j$ if $i<j$. The indegree sequence
of $T$ can be considered as a partition $\lambda \vdash n-1$. The
enumeration of trees with a given indegree sequence arises in
counting secant planes of curves in projective spaces. Recently
Ethan Cotterill conjectured a formula for the number of trees on
$[n]$ with indegree sequence corresponding to a partition
$\lambda$. In this paper we give two proofs of Cotterill's
conjecture: one is ``semi-combinatorial" based on induction, the
other is a bijective proof.



\section{Introduction}

For a labelled tree on the vertex set $[n]:=\{1,2,\ldots, n\}$,
define the direction of each edge $ij$ as $i \rightarrow j$ if
$i<j$. The indegree sequence of $T$ can be considered as a partition
$\lambda \vdash n-1$. The problem of counting labelled trees with a
given indegree sequence was encountered by Ethan Cotterill
\cite{Cotterill} when counting secant planes of curves in projective
spaces. Write $\lambda=\langle 1^{m_1}2^{m_2}\cdots\rangle$ if
$\lambda$ has $m_i$ parts equal to $i$. Given $\lambda=\langle
1^{m_1}2^{m_2}\cdots\rangle \vdash n-1$, let $k$ be the number of
parts of $\lambda$, and $a_{\lambda}$ be the number of trees on
$[n]$ with indegree sequence corresponding to $\lambda$. Cotterill
\cite[Page 29]{Cotterill} conjectured the following result:

\begin{equation}\label{Cotterill}
a_{\lambda}=\frac{(n-1)!^2}{(n-k)!{1!}^{m_1}{2!}^{m_2}\cdots{m_1}!{m_2}!\cdots}.
\end{equation}

Note that the above formula can also be written as
\begin{equation}\label{rstan}
a_{\lambda}=\frac{(n-1)!}{(n-k)!}\cdot \frac{(n-1)!}{{1!}^{m_1}{m_1}!
{2!}^{m_2}{m_2}!\cdots},
\end{equation}
in which the second factor on the right hand side counts the
number of partitions $\pi$ of an $(n-1)$-element set of {\it type}
$\lambda$, i.e., the block sizes of $\pi$ are $\lambda_1,
\lambda_2, \ldots$. This suggests that it may help to prove
(\ref{Cotterill}) if we can find a map $\phi:
\mathcal{T}_{\lambda} \rightarrow \Pi_{\lambda}$ for any $\lambda
\vdash n-1$, where $\mathcal{T}_{\lambda}$ is the set of trees on
$[n]$ with indegree sequence $\lambda$, and $\Pi_{\lambda}$ is the
set of partitions of $[2,n]:=\{2,3,\ldots, n\}$ of type $\lambda$.
Richard Stanley (personal communication) suggested that such a map
$\phi$ can be defined as follows.

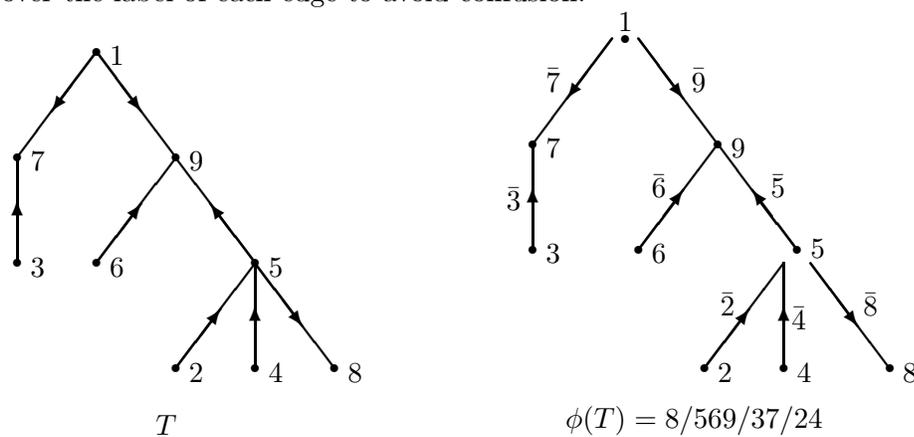
\begin{figure}[h,t]
\begin{center}
\setlength{\unitlength}{10pt}
\thicklines

\begin{picture}(35,15)
\put(8,0){\makebox(0,0)[b]{$T$}}
\put(0,2.5){
\begin{picture}(14,12)
\put(5,12){\circle*{0.3}}
\put(5.5,11.5){1}

\put(2,8){\circle*{0.3}} \put(2.5,7.5){7}
\put(5,12){\line(-3,-4){3}} \put(5,12){\vector(-3,-4){1.7}}
\put(2,4){\circle*{0.3}} \put(2.5,3.5){3}
\put(2,4){\line(0,1){4}}\put(2,4){\vector(0,1){2.4}}

\put(8,8){\circle*{0.3}} \put(8.5,7.5){9}
\put(5,12){\line(3,-4){3}} \put(5,12){\vector(3,-4){1.7}}
\put(5,4){\circle*{0.3}} \put(5.5,3.5){6}
\put(5,4){\line(3,4){3}} \put(5,4){\vector(3,4){1.7}}
\put(11,4){\circle*{0.3}} \put(11.5,3.5){5}
\put(11,4){\line(-3,4){3}}\put(11,4){\vector(-3,4){1.7}}

\put(8,0){\circle*{0.3}} \put(8.5,-0.5){2}
\put(8,0){\line(3,4){3}} \put(8,0){\vector(3,4){1.7}}
\put(11,0){\circle*{0.3}} \put(11.5,-0.5){4}
\put(11,0){\line(0,1){4}} \put(11,0){\vector(0,1){2.4}}

\put(14,0){\circle*{0.3}} \put(14.5,-0.5){8}
\put(11,4){\line(3,-4){3}} \put(11,4){\vector(3,-4){1.8}}
\end{picture}}


\put(28,0){\makebox(0,0)[b]{$\phi(T)=8/569/37/24$}}
\put(20,3){
\begin{picture}(14,12)
\put(5,12){\circle*{0.3}} \put(5,12.3){\makebox(0,0)[b]{1}}

\put(8.5,8){\circle*{0.3}} \put(9,7.5){9}
\put(5.5,4){\circle*{0.3}} \put(6,3.5){6}
\put(11.5,4){\circle*{0.3}} \put(12,3.7){5}
\put(5.5,12){\line(3,-4){3}} \put(5.5,12){\vector(3,-4){1.7}} \put(7.5,10){$\bar{9}$}
\put(5.5,4){\line(3,4){3}} \put(5.5,4){\vector(3,4){1.7}} \put(6,6){$\bar{6}$}
\put(11.5,4){\line(-3,4){3}} \put(11.5,4){\vector(-3,4){1.7}} \put(10.5,6){$\bar{5}$}

\put(1.5,8){\circle*{0.3}} \put(2,7.5){7}
\put(4.5,12){\line(-3,-4){3}} \put(4.5,12){\vector(-3,-4){1.7}} \put(2,10){$\bar{7}$}
\put(1.5,4){\circle*{0.3}} \put(2,3.5){3}
\put(1.5,4){\line(0,1){4}}\put(1.5,4){\vector(0,1){2.4}} \put(0.5,5.5){$\bar{3}$}

\put(8,-0.5){\circle*{0.3}} \put(8.5,-1){2}
\put(8,-0.5){\line(3,4){3}} \put(8,-0.5){\vector(3,4){1.7}} \put(8.6,1.5){$\bar{2}$}
\put(11,-0.5){\circle*{0.3}} \put(11.5,-1){4}
\put(11,-0.5){\line(0,1){4}} \put(11,-0.5){\vector(0,1){2.4}}\put(11.3,1){$\bar{4}$}

\put(15,-0.5){\circle*{0.3}} \put(15.5,-1){8}
\put(12,3.5){\line(3,-4){3}} \put(12,3.5){\vector(3,-4){1.8}}\put(14,1.5){$\bar{8}$}

\end{picture}}

\end{picture}
\caption{A tree $T\in \mathcal{T}_{3221}$, and $\phi(T)=8/569/37/24 \in
\Pi_{3221}$.}\label{fig_T2Pi}
\end{center}
\end{figure}

Given $\lambda \vdash n-1$ and $T \in \mathcal{T}_{\lambda}$, we
can consider $T$ as a rooted tree on $[n]$ with the root $1$
``hung up" (See Figure \ref{fig_T2Pi}). Now we label the edges of
$T$ such that each edge has the same label as the vertex right
below it. It is obvious that during the labelling each number in
$[2,n]$ is used exactly once. Putting the labels of those edges
which point to the same vertex into one block, we get a partition
$\pi \in \Pi_{\lambda}$. Figure \ref{fig_T2Pi} shows a tree $T\in
\mathcal{T}_{3221}$, and $\phi(T)=\pi=8/569/37/24 \in \Pi_{3221}$.
We put a bar over the label of each edge to avoid confusion.

While the map $\phi$ gives a natural interpretation of the second
factor in equation (\ref{rstan}), one can easily check that the
preimage of $\phi$ is not unique: we can get the same partition by
applying $\phi$ to different trees. Let $\mathcal{T}_\pi$ be the
set of preimages of $\pi\in\Pi_\lambda$ under the map $\phi$,
i.e., $\mathcal{T}_\pi=\phi^{-1}(\pi)$, and let
$f(\pi):=|\mathcal{T}_\pi|$. Then
$\mathcal{T}_{\lambda}=\bigcup_{\pi \in
\Pi_{\lambda}}\mathcal{T}_{\pi}$. Our main task is to prove the
following theorem.

\begin{theo}\label{main_thm}
Given $\lambda \vdash n-1$ and $\pi \in \Pi_{\lambda}$, we
have
\begin{equation*}
f(\pi)=|\mathcal{T}_\pi|=\frac{(n-1)!}{(n-|\pi|)!},
\end{equation*}
where $|\pi|$ is the number of blocks of $\pi$.
\end{theo}

In the remainder of this paper we give proofs of this result using
two different approaches. In Section 2, we give a
``semi-combinatorial" proof based on induction on $n$. In Section
3 and Section 4, we give a bijective proof. Finally in section 5,
some further problems are raised.

\section{A Semi-combinatorial Proof}
In this section, we will give an inductive proof of Theorem
\ref{main_thm}.

\begin{lem}\label{lemma_for_induction}
The value $f(\pi)$ is independent of $\pi\in\Pi_\lambda$, i.e.,
for any $\pi_1,\pi_2\in\Pi_\lambda$, we have
$f({\pi_1})=f({\pi_2})$.
\end{lem}

\pf Since the symmetric group of $[2,n]$ is generated by {\it
adjacent transpositions} $\{s_i:2\leq i\leq n-1\}$, where
$s_i=(i,i+1)$ is the function that swaps two elements $i$ and
$i+1$, it suffices to show that $f({\pi_1})=f({\pi_2})$ for any
$\pi_1,\pi_2\in\Pi_\lambda$ such that by switching $i$ and $i+1$
in $\pi_2$ we will get $\pi_1$ ($2\leq i\leq n-1$). If $i$ and
$i+1$ are in the same block of $\pi_1$, then $\pi_1=\pi_2$. The
assertion is trivial in this case. In the following, we will
assume that $i$ and $i+1$ are in different blocks of $\pi_1$.

In order to prove $f({\pi_1})=f({\pi_2})$, we construct an
involution
$\varphi_i:\mathcal{T}_{\pi_1}\cup\mathcal{T}_{\pi_2}\rightarrow\mathcal{T}_{\pi_1}\cup\mathcal{T}_{\pi_2}$.
For any tree $T\in\mathcal{T}_{\pi_1}\cup\mathcal{T}_{\pi_2}$,
consider the two vertices labelled $i$ and $i+1$.

If vertices $i$ and $i+1$ are not adjacent, exchanging the labels
of these two vertices will give us a new tree $T'$. Let
$\varphi_i(T)=T'$.

\begin{figure}[h,t]
\begin{center}
\setlength{\unitlength}{10pt}
\begin{picture}(17,14)
\put(0,0){\line(0,1){1}}\put(2,0.5){\vector(-1,0){2}}
\put(6,0){\line(0,1){1}}\put(4,0.5){\vector(1,0){2}}
\put(3,0.5){\makebox(0,0){$T_{i}$}}

\put(11,0){\line(0,1){1}}\put(12.8,0.5){\vector(-1,0){1.8}}
\put(17,0){\line(0,1){1}}\put(15.2,0.5){\vector(1,0){1.8}}
\put(14,0.5){\makebox(0,0){$T_{i+1}$}}

\put(0,2){
\begin{picture}(17,12)
\put(2,2){\circle{4}} \put(2,2){\makebox(0,0){$B_i$}}
\put(15,2){\circle{4}} \put(15,2){\makebox(0,0){$B_{i+1}$}}
\put(2,10){\circle{4}} \put(2,10){\makebox(0,0){$A_i$}}
\put(15,10){\circle{4}} \put(15,10){\makebox(0,0){$A_{i+1}$}}
\thicklines \put(6,6){\line(1,0){5}}\put(6,6){\vector(1,0){2.7}}
\put(6,6){\circle*{0.3}}\put(11,6){\circle*{0.3}}
\put(5.4,6){\makebox(0,0)[r]{$i$}}
\put(11.8,6){\makebox(0,0)[l]{$i\!+\!\!1$}}

\thinlines
\put(6,6){\line(-3,2){3.2}}\put(6,6){\vector(-3,2){2}}
\put(6,6){\line(-1,1){2.56}}\put(6,6){\vector(-1,1){1.7}}
\put(6,6){\line(-2,3){2.15}}\put(6,6){\vector(-2,3){1.35}}

\put(11,6){\line(3,2){3.2}}\put(11,6){\vector(3,2){2}}
\put(11,6){\line(1,1){2.56}}\put(11,6){\vector(1,1){1.7}}
\put(11,6){\line(2,3){2.15}}\put(11,6){\vector(2,3){1.35}}

\put(6,6){\line(-3,-2){3.2}}\put(2.8,3.86){\vector(3,2){1.6}}
\put(6,6){\line(-1,-1){2.56}} \put(3.44,3.44){\vector(1,1){1.28}}
\put(6,6){\line(-2,-3){2.15}}\put(3.85,2.78){\vector(2,3){1.08}}

\put(11,6){\line(3,-2){3.2}}   \put(14.2,3.86){\vector(-3,2){1.6}}
\put(11,6){\line(1,-1){2.56}}  \put(13.56,3.44){\vector(-1,1){1.28}}
\put(11,6){\line(2,-3){2.15}}  \put(13.15,2.78){\vector(-2,3){1.08}}
\end{picture}}
\end{picture}
\caption{A partition of the tree $T$.} \label{fig:partition}
\end{center}
\end{figure}
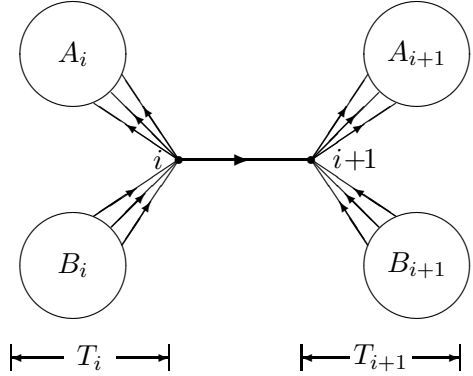

If vertices $i$ and $i+1$ are adjacent, let $T_i$ (resp.
$T_{i+1}$) be the largest subtree containing vertex $i$ but not
$i+1$ (resp. containing $i+1$ but not $i$). For $j=i,i+1$, let
$T_j=\{j\}\cup A_j\cup B_j$, where $A_j$ (resp. $B_j$) is the
sub-forest such that every edge between itself and vertex $j$ is
pointing away from $j$ (resp. pointing to $j$). (See Figure
\ref{fig:partition}.)

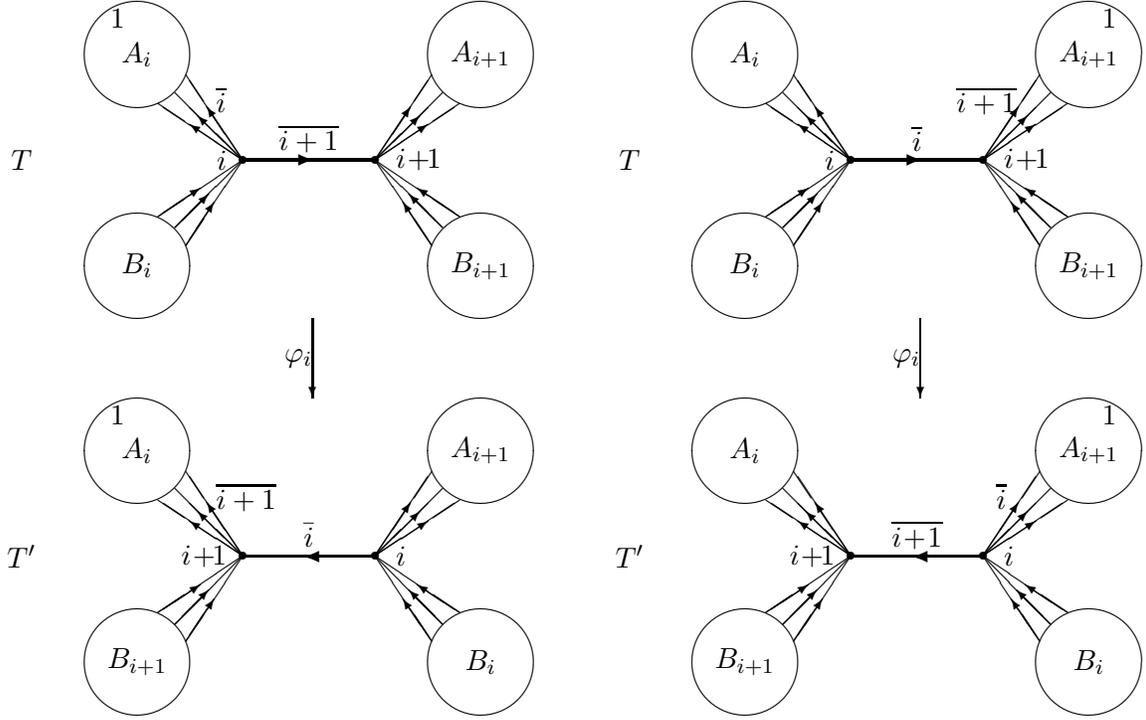
\begin{figure}[h,t]
\begin{center}
\setlength{\unitlength}{10pt}
\begin{picture}(40,27)

\put(0,6){\makebox(0,0){$T^{\prime}$}}
\put(2,0){
\begin{picture}(17,12)
\put(2,2){\circle{4}} \put(2,2){\makebox(0,0){$B_{i+1}$}}
\put(15,2){\circle{4}} \put(15,2){\makebox(0,0){$B_{i}$}}
\put(2,10){\circle{4}} \put(2,10){\makebox(0,0){$A_i$}}
\put(1,11){1} \put(15,10){\circle{4}}
\put(15,10){\makebox(0,0){$A_{i+1}$}} \thicklines
\put(6,6){\line(1,0){5}}\put(11,6){\vector(-1,0){2.7}}
\put(6,6){\circle*{0.3}}\put(11,6){\circle*{0.3}}
\put(5.4,6){\makebox(0,0)[r]{$i\!+\!\!1$}}
\put(11.8,6){\makebox(0,0)[l]{$i$}}
\put(8.5,6.8){\makebox(0,0){$\bar{i}$}} \thinlines
\put(5,7.8){$\overline{i+1}$}

\put(6,6){\line(-3,2){3.2}}\put(6,6){\vector(-3,2){2}}
\put(6,6){\line(-1,1){2.56}}\put(6,6){\vector(-1,1){1.7}}
\put(6,6){\line(-2,3){2.15}}\put(6,6){\vector(-2,3){1.35}}
\put(11,6){\line(3,2){3.2}}\put(11,6){\vector(3,2){2}}
\put(11,6){\line(1,1){2.56}}\put(11,6){\vector(1,1){1.7}}
\put(11,6){\line(2,3){2.15}}\put(11,6){\vector(2,3){1.35}}

\put(6,6){\line(-3,-2){3.2}}\put(2.8,3.86){\vector(3,2){1.6}}
\put(6,6){\line(-1,-1){2.56}} \put(3.44,3.44){\vector(1,1){1.28}}
\put(6,6){\line(-2,-3){2.15}}\put(3.85,2.78){\vector(2,3){1.08}}

\put(11,6){\line(3,-2){3.2}}   \put(14.2,3.86){\vector(-3,2){1.6}}
\put(11,6){\line(1,-1){2.56}}  \put(13.56,3.44){\vector(-1,1){1.28}}
\put(11,6){\line(2,-3){2.15}}  \put(13.15,2.78){\vector(-2,3){1.08}}
\end{picture}}

\put(11,15){\vector(0,-1){3}}
\put(11,13.5){\makebox(0,0)[r]{$\varphi_i$}}

\put(0,21){\makebox(0,0){$T$}}
\put(2,15){
\begin{picture}(17,12)
\put(2,2){\circle{4}} \put(2,2){\makebox(0,0){$B_i$}}
\put(15,2){\circle{4}} \put(15,2){\makebox(0,0){$B_{i+1}$}}
\put(2,10){\circle{4}} \put(2,10){\makebox(0,0){$A_i$}}
\put(1,11){1} \put(15,10){\circle{4}}
\put(15,10){\makebox(0,0){$A_{i+1}$}} \thicklines
\put(6,6){\line(1,0){5}}\put(6,6){\vector(1,0){2.7}}
\put(6,6){\circle*{0.3}}\put(11,6){\circle*{0.3}}
\put(5.4,6){\makebox(0,0)[r]{$i$}}
\put(11.8,6){\makebox(0,0)[l]{$i\!+\!\!1$}}
\put(8.5,6.8){\makebox(0,0){$\overline{i+1}$}} \thinlines

\put(6,6){\line(-3,2){3.2}}\put(6,6){\vector(-3,2){2}}
\put(6,6){\line(-1,1){2.56}}\put(6,6){\vector(-1,1){1.7}}
\put(6,6){\line(-2,3){2.15}}\put(6,6){\vector(-2,3){1.35}}
\put(11,6){\line(3,2){3.2}}\put(11,6){\vector(3,2){2}}
\put(11,6){\line(1,1){2.56}}\put(11,6){\vector(1,1){1.7}}
\put(11,6){\line(2,3){2.15}}\put(11,6){\vector(2,3){1.35}}

\put(5,7.8){$\overline{i}$}

\put(6,6){\line(-3,-2){3.2}}\put(2.8,3.86){\vector(3,2){1.6}}
\put(6,6){\line(-1,-1){2.56}} \put(3.44,3.44){\vector(1,1){1.28}}
\put(6,6){\line(-2,-3){2.15}}\put(3.85,2.78){\vector(2,3){1.08}}

\put(11,6){\line(3,-2){3.2}}   \put(14.2,3.86){\vector(-3,2){1.6}}
\put(11,6){\line(1,-1){2.56}}  \put(13.56,3.44){\vector(-1,1){1.28}}
\put(11,6){\line(2,-3){2.15}}  \put(13.15,2.78){\vector(-2,3){1.08}}
\end{picture}}

\put(23,6){\makebox(0,0){$T^{\prime}$}}
\put(25,0){
\begin{picture}(17,12)
\put(2,2){\circle{4}} \put(2,2){\makebox(0,0){$B_{i+1}$}}
\put(15,2){\circle{4}} \put(15,2){\makebox(0,0){$B_{i}$}}
\put(2,10){\circle{4}} \put(2,10){\makebox(0,0){$A_i$}}
\put(15,10){\circle{4}}
\put(15,10){\makebox(0,0){$A_{i+1}$}}\put(15.5,11){1} \thicklines
\put(6,6){\line(1,0){5}}\put(11,6){\vector(-1,0){2.7}}
\put(6,6){\circle*{0.3}}\put(11,6){\circle*{0.3}}
\put(5.4,6){\makebox(0,0)[r]{$i\!+\!\!1$}}
\put(11.8,6){\makebox(0,0)[l]{$i$}}
\put(8.5,6.8){\makebox(0,0){$\overline{i\!+\!1}$}} \thinlines
\put(11.5,7.8){$\overline{i}$}
\put(6,6){\line(-3,2){3.2}}\put(6,6){\vector(-3,2){2}}
\put(6,6){\line(-1,1){2.56}}\put(6,6){\vector(-1,1){1.7}}
\put(6,6){\line(-2,3){2.15}}\put(6,6){\vector(-2,3){1.35}}
\put(11,6){\line(3,2){3.2}}\put(11,6){\vector(3,2){2}}
\put(11,6){\line(1,1){2.56}}\put(11,6){\vector(1,1){1.7}}
\put(11,6){\line(2,3){2.15}}\put(11,6){\vector(2,3){1.35}}

\put(6,6){\line(-3,-2){3.2}}\put(2.8,3.86){\vector(3,2){1.6}}
\put(6,6){\line(-1,-1){2.56}} \put(3.44,3.44){\vector(1,1){1.28}}
\put(6,6){\line(-2,-3){2.15}}\put(3.85,2.78){\vector(2,3){1.08}}

\put(11,6){\line(3,-2){3.2}}   \put(14.2,3.86){\vector(-3,2){1.6}}
\put(11,6){\line(1,-1){2.56}}  \put(13.56,3.44){\vector(-1,1){1.28}}
\put(11,6){\line(2,-3){2.15}}  \put(13.15,2.78){\vector(-2,3){1.08}}
\end{picture}}

\put(34,15){\vector(0,-1){3}}
\put(34,13.5){\makebox(0,0)[r]{$\varphi_i$}}

\put(23,21){\makebox(0,0){$T$}}
\put(25,15){
\begin{picture}(17,12)
\put(2,2){\circle{4}} \put(2,2){\makebox(0,0){$B_i$}}
\put(15,2){\circle{4}} \put(15,2){\makebox(0,0){$B_{i+1}$}}
\put(2,10){\circle{4}} \put(2,10){\makebox(0,0){$A_i$}}
\put(15,10){\circle{4}}
\put(15,10){\makebox(0,0){$A_{i+1}$}}\put(15.5,11){1} \thicklines
\put(6,6){\line(1,0){5}}\put(6,6){\vector(1,0){2.7}}
\put(6,6){\circle*{0.3}}\put(11,6){\circle*{0.3}}
\put(5.4,6){\makebox(0,0)[r]{$i$}}
\put(11.8,6){\makebox(0,0)[l]{$i\!+\!\!1$}}
\put(8.5,6.8){\makebox(0,0){$\overline{i}$}} \thinlines
\put(10,7.8){$\overline{i+1}$}

\put(6,6){\line(-3,2){3.2}}\put(6,6){\vector(-3,2){2}}
\put(6,6){\line(-1,1){2.56}}\put(6,6){\vector(-1,1){1.7}}
\put(6,6){\line(-2,3){2.15}}\put(6,6){\vector(-2,3){1.35}}
\put(11,6){\line(3,2){3.2}}\put(11,6){\vector(3,2){2}}
\put(11,6){\line(1,1){2.56}}\put(11,6){\vector(1,1){1.7}}
\put(11,6){\line(2,3){2.15}}\put(11,6){\vector(2,3){1.35}}

\put(6,6){\line(-3,-2){3.2}}\put(2.8,3.86){\vector(3,2){1.6}}
\put(6,6){\line(-1,-1){2.56}} \put(3.44,3.44){\vector(1,1){1.28}}
\put(6,6){\line(-2,-3){2.15}}\put(3.85,2.78){\vector(2,3){1.08}}

\put(11,6){\line(3,-2){3.2}}   \put(14.2,3.86){\vector(-3,2){1.6}}
\put(11,6){\line(1,-1){2.56}}  \put(13.56,3.44){\vector(-1,1){1.28}}
\put(11,6){\line(2,-3){2.15}}  \put(13.15,2.78){\vector(-2,3){1.08}}
\end{picture}}
\end{picture}
\caption{Map $\varphi_i$ (left: $1$ in $A_i$, right: $1$ in
$A_{i+1}$).} \label{fig:1inA}
\end{center}
\end{figure}

Considering the position of vertex $1$, there are three cases:

\noindent{\it Case 1:} If vertex $1$ is in either $A_i$ or
$A_{i+1}$, make all edges from $B_i$ to vertex $i$ point to vertex
$i+1$ instead, make all edges from $B_{i+1}$ to vertex $i+1$ point
to vertex $i$ instead, and switch the vertex labels $i$ and $i+1$
(at the same time the direction of the edge between $i$ and $i+1$
will be changed automatically). Then we will get a new tree $T'$.
Let $\varphi_i(T)=T'$. (See Figure \ref{fig:1inA})

\begin{figure}[htbp]
\begin{center}
\setlength{\unitlength}{10pt}
\begin{picture}(40,27)

\put(11,0){\makebox(0,0)[]{(1)}}
\put(0,8){\makebox(0,0){$T^{\prime}$}}
\put(2,2){
\begin{picture}(17,12)
\put(11,1.5){\circle{3}} \put(11,1.8){\makebox(0,0){$B_{i}^{\prime}$}}
\put(11,0.1){\makebox(0,0)[b]{\small 1}}
\put(11,3.1){\line(0,1){3}}
\put(11,3.1){\vector(0,1){1.7}}
\put(10.5,4){$\bar{i}$}

\put(2,2){\circle{4}} \put(2,2){\makebox(0,0){$B_{i}^{\prime\prime}$}}
\put(15,2){\circle{4}} \put(15,2){\makebox(0,0){$B_{i+1}$}}
\put(2,10){\circle{4}} \put(2,10){\makebox(0,0){$A_i$}}
\put(15,10){\circle{4}} \put(15,10){\makebox(0,0){$A_{i+1}$}}
\thicklines
\put(6,6){\line(1,0){5}}\put(11,6){\vector(-1,0){2.7}}
\put(6,6){\circle*{0.3}}\put(11,6){\circle*{0.3}}
\put(5.3,6){\makebox(0,0)[r]{$i\!+\!\!1$}}
\put(11.8,6){\makebox(0,0)[l]{$i$}}
\put(8.5,6.8){\makebox(0,0){$\overline{i+1}$}}
\thinlines
\put(6,6){\line(-3,2){3.2}}\put(6,6){\vector(-3,2){2}}
\put(6,6){\line(-1,1){2.56}}\put(6,6){\vector(-1,1){1.7}}
\put(6,6){\line(-2,3){2.15}}\put(6,6){\vector(-2,3){1.35}}
\put(11,6){\line(3,2){3.2}}\put(11,6){\vector(3,2){2}}
\put(11,6){\line(1,1){2.56}}\put(11,6){\vector(1,1){1.7}}
\put(11,6){\line(2,3){2.15}}\put(11,6){\vector(2,3){1.35}}

\put(6,6){\line(-3,-2){3.2}}\put(2.8,3.86){\vector(3,2){1.6}}
\put(6,6){\line(-1,-1){2.56}} \put(3.44,3.44){\vector(1,1){1.28}}

\put(11,6){\line(3,-2){3.2}}   \put(14.2,3.86){\vector(-3,2){1.6}}
\put(11,6){\line(1,-1){2.56}}  \put(13.56,3.44){\vector(-1,1){1.28}}
\put(11,6){\line(2,-3){2.15}}  \put(13.15,2.78){\vector(-2,3){1.08}}
\end{picture}}

\put(11,17){\vector(0,-1){3}}
\put(11,15.5){\makebox(0,0)[r]{$\varphi_i$}}

\put(0,23){\makebox(0,0){$T$}}
\put(2,17){
\begin{picture}(17,12)
\put(6,1.5){\circle{3}} \put(6,1.8){\makebox(0,0){$B_{i}^{\prime}$}}
\put(6,0.1){\makebox(0,0)[b]{\small 1}}
\put(6,3.1){\line(0,1){3}}
\put(6,3.1){\vector(0,1){1.7}}
\put(6.4,4){$\bar{i}$}

\put(2,2){\circle{4}} \put(2,2){\makebox(0,0){$B_{i}^{\prime\prime}$}}
\put(15,2){\circle{4}} \put(15,2){\makebox(0,0){$B_{i+1}$}}
\put(2,10){\circle{4}} \put(2,10){\makebox(0,0){$A_i$}}
\put(15,10){\circle{4}} \put(15,10){\makebox(0,0){$A_{i+1}$}}
\thicklines
\put(6,6){\line(1,0){5}}\put(6,6){\vector(1,0){2.7}}
\put(6,6){\circle*{0.3}}\put(11,6){\circle*{0.3}}
\put(5.2,6){\makebox(0,0)[r]{$i$}}
\put(11.8,6){\makebox(0,0)[l]{$i\!+\!\!1$}}
\put(8.5,6.8){\makebox(0,0){$\overline{i+1}$}}
\thinlines
\put(6,6){\line(-3,2){3.2}}\put(6,6){\vector(-3,2){2}}
\put(6,6){\line(-1,1){2.56}}\put(6,6){\vector(-1,1){1.7}}
\put(6,6){\line(-2,3){2.15}}\put(6,6){\vector(-2,3){1.35}}
\put(11,6){\line(3,2){3.2}}\put(11,6){\vector(3,2){2}}
\put(11,6){\line(1,1){2.56}}\put(11,6){\vector(1,1){1.7}}
\put(11,6){\line(2,3){2.15}}\put(11,6){\vector(2,3){1.35}}

\put(6,6){\line(-3,-2){3.2}}\put(2.8,3.86){\vector(3,2){1.6}}
\put(6,6){\line(-1,-1){2.56}} \put(3.44,3.44){\vector(1,1){1.28}}

\put(11,6){\line(3,-2){3.2}}   \put(14.2,3.86){\vector(-3,2){1.6}}
\put(11,6){\line(1,-1){2.56}}  \put(13.56,3.44){\vector(-1,1){1.28}}
\put(11,6){\line(2,-3){2.15}}  \put(13.15,2.78){\vector(-2,3){1.08}}
\end{picture}}

\put(34,8){\makebox(0,0)[]{(2)}}
\put(23,15.5){\makebox(0,0){$T$}}
\put(25,9.5){
\begin{picture}(17,12)
\put(2,2){\circle{4}} \put(2,2){\makebox(0,0){$B_{i}$}}
\put(15,2){\circle{4}} \put(15,2){\makebox(0,0){$B_{i+1}$}}\put(15,0.2){\small 1}
\put(2,10){\circle{4}} \put(2,10){\makebox(0,0){$A_i$}}
\put(15,10){\circle{4}} \put(15,10){\makebox(0,0){$A_{i+1}$}}
\thicklines
\put(6,6){\line(1,0){5}}\put(6,6){\vector(1,0){2.7}}
\put(6,6){\circle*{0.3}}\put(11,6){\circle*{0.3}}
\put(5.2,6){\makebox(0,0)[r]{$i$}}
\put(11.8,6){\makebox(0,0)[l]{$i\!+\!\!1$}}
\put(8.5,6.8){\makebox(0,0){$\bar{i}$}}
\thinlines
\put(6,6){\line(-3,2){3.2}}\put(6,6){\vector(-3,2){2}}
\put(6,6){\line(-1,1){2.56}}\put(6,6){\vector(-1,1){1.7}}
\put(6,6){\line(-2,3){2.15}}\put(6,6){\vector(-2,3){1.35}}
\put(11,6){\line(3,2){3.2}}\put(11,6){\vector(3,2){2}}
\put(11,6){\line(1,1){2.56}}\put(11,6){\vector(1,1){1.7}}
\put(11,6){\line(2,3){2.15}}\put(11,6){\vector(2,3){1.35}}

\put(6,6){\line(-3,-2){3.2}}\put(2.8,3.86){\vector(3,2){1.6}}
\put(6,6){\line(-1,-1){2.56}} \put(3.44,3.44){\vector(1,1){1.28}}
\put(6,6){\line(-2,-3){2.15}}\put(3.85,2.78){\vector(2,3){1.08}}

\put(10,3.4){$\overline{i\!+\!1}$}
\put(11,6){\line(3,-2){3.2}}   \put(14.2,3.86){\vector(-3,2){1.6}}
\put(11,6){\line(1,-1){2.56}}  \put(13.56,3.44){\vector(-1,1){1.28}}
\put(11,6){\line(2,-3){2.15}}  \put(13.15,2.78){\vector(-2,3){1.08}}
\end{picture}}

\end{picture}
\caption{Map $\varphi_i$ (left: $1$ in $B_i$, right: $1$ in
$B_{i+1}$(impossible)).} \label{fig:1inB}
\end{center}
\end{figure}
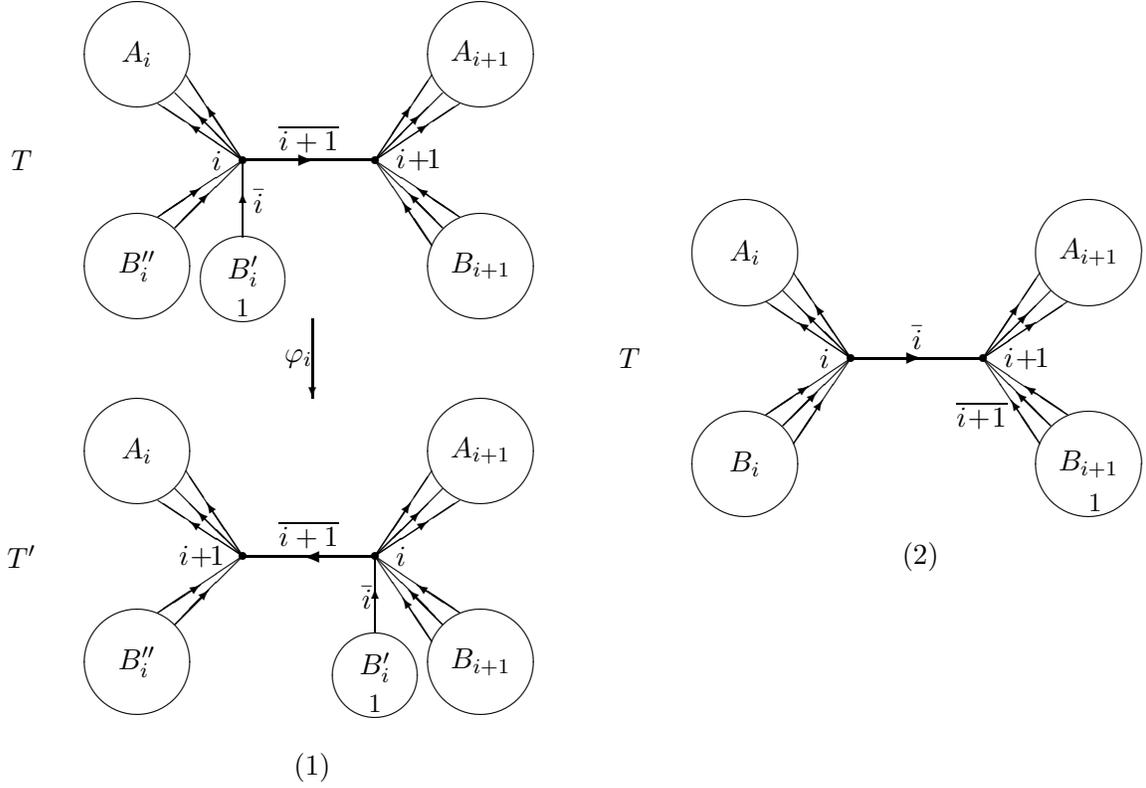

\noindent{\it Case 2:} If vertex $1$ is in $B_i$, let $B'_i$ be
the maximum subtree of $B_i$ which contains vertex $1$, and let
$B''_i$ be $B_i\backslash B_i'$. Make all edges from $B'_i$ to
vertex $i$ point to vertex $i+1$ instead, and switch the vertex
labels $i$ and $i+1$ (at the same time the direction of the edge
between $i$ and $i+1$ will be changed automatically). Then we will
get a new tree $T'$. Let $\varphi_i(T)=T'$. (See Figure
\ref{fig:1inB}(1)).

\noindent{\it Case 3:} If vertex $1$ is in $B_{i+1}$, both edges
labelled $\overline{i}$ and $\overline{i+1}$ are pointing to vertex
$i+1$, i.e., $i$ and $i+1$ are in the same block of $\pi_1$ or
$\pi_2$, then we have a contradiction to the assumption. (See Figure
\ref{fig:1inB}(2)).

From the definition of the map, we can easily check that $\phi(T)$
and $\phi(T')$ only differ in the positions of $i$ and $i+1$,
i.e., $\phi(T')$ is the same as $\phi(T)$ after switching $i$ and
$i+1$. Since $\varphi_i(T)=T'$, we have
$\varphi_i:\mathcal{T}_{\pi_1}\cup\mathcal{T}_{\pi_2}\rightarrow\mathcal{T}_{\pi_1}\cup\mathcal{T}_{\pi_2}$
is well-defined, and
$\varphi_i(\mathcal{T}_{\pi_1})\in\mathcal{T}_{\pi_2},\varphi_i(\mathcal{T}_{\pi_2})\in\mathcal{T}_{\pi_1}$.
And by applying $\varphi_i$ again, we have
$\varphi_i(\varphi_i(T))=T$. Hence, $\varphi_i$ is an involution
with no fixed points. Hence, we have
$|\mathcal{T}_{\pi_1}|=|\mathcal{T}_{\pi_2}|$, i.e.,
$f({\pi_1})=f({\pi_2})$.\qed

\noindent \emph{Proof of Theorem \ref{main_thm}: }Now with Lemma
\ref{lemma_for_induction} we can prove Theorem \ref{main_thm} by
induction on $n$, the number of vertices.

Let $\lambda=(\lambda_1,\lambda_2,\ldots,\lambda_k) \vdash n-1$,
where $\lambda_1\geq\lambda_2\geq\cdots \geq\lambda_k\geq 1$. Then
what we need to show is that for any $\pi \in \Pi_{\lambda}$, we
have $f(\pi)=(n-1)!/(n-k)!$.

\noindent \emph{Base case: }If $n=1$, we have
$k=0,\lambda=\emptyset,\pi=\emptyset$, and
$f(\pi)=1=(n-1)!/(n-k)!$.

\noindent \emph{Inductive Step: }Assume that the theorem is true
for $n-1$ $(\geq 1)$. Then consider the case for $n$.

If $\lambda_1=1$, then $\lambda=\langle 1^{n-1}\rangle$,
$\pi=n/n-1/\ldots/2$ and $k=n-1$. In this case, each
$T\in\mathcal{T}_{\pi}$ is an increasing tree, i.e., the label of
any vertex is bigger than the label of its parent, i.e., the
directions of edges are pointing away from the root $1$. Otherwise,
there is at least one vertex with indegree at least $2$,
contradicting that $\lambda$ is the indegree sequence. Hence, we can
do the bijection as in \cite[\S 1.3]{EC1} by mapping $T$ to a
permutation of $[2,n]$. Or we can use the bijection between labelled
trees and Pr\"{u}fer codes (see, for example, \cite[\S 2.4]{bondy},
or a more generalized forest version \cite[\S 5.3]{EC2}). But while
doing this bijection, what we will get is a subset of all possible
Pr\"{u}fer codes, i.e., a subset of
$[n-1]\times[n-2]\times\cdots\times[2]\times[1]$. Both methods show
that $f(\pi)=(n-1)!=(n-1)!/(n-k)!$.

Now suppose that $\lambda_1\geq 2$. By Lemma
\ref{lemma_for_induction}, we can assume without loss of generality,
that both $n$ and $n-1$ are in the same block $B_1$ of
$\pi=\{B_1,B_2,\ldots,B_k\}$. Pick $T\in\mathcal{T}_{\pi}$. Since
$n$ is the largest label, by the definition of $\pi$ and
$\mathcal{T}_{\pi}$, we know that vertices $n$ and $n-1$ are
adjacent. By merging the edge between $n$ and $n-1$ in $T$, and
deleting the label $n$, we get a new tree $\tilde{T}$ with $n-1$
vertices. There are two possible cases:

\noindent{\it Case 1:} If the indegree of vertex $n-1$ in $T$ is
$0$, then
$\phi(\tilde{T})=\{B_1\backslash\{n\},B_2,\ldots,B_k\}=:\tilde{\pi}_1$.

\noindent{\it Case 2:} If the indegree of vertex $n-1$ in $T$ is
not $0$, then there exists $j\in [2,k]$ such that
$\phi(\tilde{T})=\{B_1\cup
B_j\backslash\{n\},B_2,\ldots,B_{j-1},B_{j+1},\ldots,B_k\}=:\tilde{\pi}_j$.

One can easily check that this is a bijection. Thus,
$f(\pi)=\sum_{j=1}^{k}f({\tilde{\pi}_j})$. By the induction
hypothesis we have
\[f(\pi)=\sum_{j=1}^{k}f({\tilde{\pi}_j})=\frac{((n-1)-1)!}{((n-1)-k)!}+(k-1)\frac{((n-1)-1)!}{((n-1)-(k-1))!}=\frac{(n-1)!}{(n-k)!},\]
which proves the case for $n$.

Hence it follows by induction that Theorem \ref{main_thm} is true
for all possible $n$.\qed

\section{An ``Almost" Bijective Proof}
The inductive proof in the former section makes Corterill's
conjecture a theorem, but it does not explain combinatorially why
there is such a simple factor $(n-1)!/(n-k)!$. In this section, we
will try to give a bijective proof to explain this fact.

First we will give some terminology and notation related to
posets. Let $S$ be a finite set. We use $\Pi_{S}$ to denote the
poset (actually a geometric lattice) of all partitions of $S$
ordered by refinement ($\sigma \preceq \pi$ in $\Pi_{S}$ if every
block of $\sigma$ is contained in a block of $\pi$). In the
following discussion we will consider the case that $S=[2,n]$.

Second, we will state the basic definitions. Given
$\pi\in\Pi_{[2,n]}$, recall that $\mathcal{T}_\pi$ is the set of
labelled trees with preimage $\pi$ under the map $\phi$. Let $B,
B'$ be two subsets of $[2,n]$. We say that $B\leq B'$( resp.
$B<B'$) if and only if $\min B\leq\min B'$ (resp. $\min B<\min
B'$). Given $T\in\mathcal{T}_\pi$ and $\pi=\phi(T)$, let
$B=\{b_1,b_2,\ldots,b_t\}_<$ be a subset of one of the blocks of
$\pi$. We define the \emph{Star corresponding to $B$} to be the
subset of $T$ that contains all vertices and edges with labels in
the set $B$, and denote it as $\Star(B)$. Induced by the ordering
of the subsets of $[2,n]$, we will also get an ordering of the
stars. For $\Star(B)$, there exists a unique vertex of $T$ with
some label, say $c$, such that the vertex $c$ is attached to one
of the edges in $\Star(B)$, but $c\not\in B$. We call the vertex
$c$ the \emph{cut point} of $B$, and denote it by $c(B)$.

For $T\in\mathcal{T}_\pi$ and $\sigma\preceq\pi$, we define the
{\it decomposition of $T$ with respect to
$\sigma=\{B_1,B_2,\ldots,B_k\}$} to be
$T=(\bigcup_{j=1}^{k}\Star(B_j))\cup\{$vertex $1\}$, where
$\Star(B_j)$ are the stars corresponding to $B_j$ in $T$. In this
decomposition, the \emph{leaf-stars} are the stars that don't
contain any cut points, i.e., if you remove a leaf-star from $T$,
what's left is still a connected tree.

For example, for the tree $T$ in Figure \ref{fig_T2Pi} we have
$\phi(T)=\pi=8/569/37/24$. $\Star(\{3,7\})$, $\Star(\{2,4\})$ and
$\Star(\{8\})$ are all leaf-stars of $T$, and we have
$c(\{3,7\})=1$, $c(\{2,4\})=5$, $c(\{8\})=5$ and $c(\{5,6,9\})=1$.

Now we define a variant of the map $\phi$, which turns out to be a
bijection. For any $\sigma=\{B_1,B_2,\ldots,B_k\}\in \Pi_{[2,n]}$,
let
$\mathcal{T}_{\succeq\sigma}=\bigcup_{\pi\succeq\sigma}\mathcal{T}_\pi$.
We define $\phi_\sigma:\mathcal{T}_{\succeq\sigma}\rightarrow
[n]^{k-1}$ as follows.
\begin{enumerate}

\item{Let $T_0=T$.}

\item{For $i=1,2,\ldots,k$, let $\Star(B(i))$ be the largest
leaf-star in the decomposition of $T_{i-1}$ with respect to
$\sigma\backslash\{B(1),B(2),\ldots,B(i-1)\}$. Then we remove
$\Star(B(i))$ and keep a record of the vertex it is attached to,
i.e., let $\omega_i=c(B(i)), T_i=T_{i-1}\backslash\Star(B(i))$.}

\end{enumerate}

\noindent Let
$\phi_\sigma(T)=\omega:=\omega_1\omega_2\ldots\omega_{k-1}\in
[n]^{k-1}$. ( We don't need to include $\omega_k$ since it is
always $1$.)

\begin{theo}
For any $\sigma\in \Pi_{[2,n]}$, the map $\phi_\sigma$ is a
bijection between $\mathcal{T}_{\succeq\sigma}$ and
$[n]^{|\sigma|-1}$.
\end{theo}

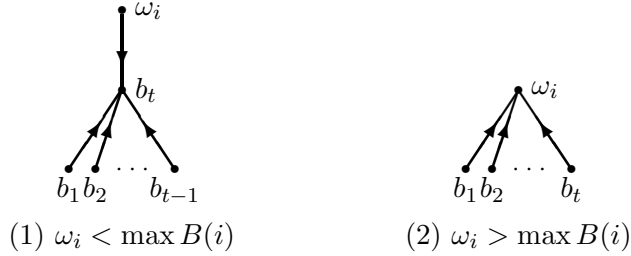
\begin{figure}[h,t]
\begin{center}
\setlength{\unitlength}{10pt}
\begin{picture}(20,10)
\thicklines

\put(0,1){
\begin{picture}(5,8)
\put(3,0){\makebox(0,0)[t]{(1) $\omega_i<\max B(i)$}}

\put(3,5){\circle*{0.3}} \put(3.5,5){\makebox(0,0)[l]{$b_t$}}

\put(1,2){\circle*{0.3}} \put(1,1.7){\makebox(0,0)[t]{$b_1$}}
\put(1,2){\line(2,3){2}} \put(1,2){\vector(2,3){1.2}}

\put(2,2){\circle*{0.3}} \put(2,1.7){\makebox(0,0)[t]{$b_2$}}
\put(2,2){\line(1,3){1}} \put(2,2){\vector(1,3){0.6}}

\multiput(2.9,2)(0.5,0){3}{\circle*{0.1}}

\put(5,2){\circle*{0.3}} \put(5,1.7){\makebox(0,0)[t]{$b_{t-1}$}}
\put(5,2){\line(-2,3){2}} \put(5,2){\vector(-2,3){1.2}}

\put(3,8){\circle*{0.3}} \put(3.5,8){\makebox(0,0)[l]{$\omega_i$}}
\put(3,7.85){\line(0,-1){3}} \put(3,7.85){\vector(0,-1){2}}
\end{picture}}

\put(15,1){
\begin{picture}(5,8)
\put(3,0){\makebox(0,0)[t]{(2) $\omega_i>\max B(i)$}}
\put(3,5){\circle*{0.3}} \put(3.5,5){\makebox(0,0)[l]{$\omega_i$}}

\put(1,2){\circle*{0.3}} \put(1,1.7){\makebox(0,0)[t]{$b_1$}}
\put(1,2){\line(2,3){1.92}} \put(1,2){\vector(2,3){1.2}}

\put(2,2){\circle*{0.3}} \put(2,1.7){\makebox(0,0)[t]{$b_2$}}
\put(2,2){\line(1,3){0.94}} \put(2,2){\vector(1,3){0.6}}

\multiput(2.9,2)(0.5,0){3}{\circle*{0.1}}

\put(5,2){\circle*{0.3}} \put(5,1.7){\makebox(0,0)[t]{$b_{t}$}}
\put(5,2){\line(-2,3){1.92}} \put(5,2){\vector(-2,3){1.2}}
\end{picture}}
\end{picture}
\caption{Two cases when attaching $B(i)$ to
$\omega_i$.}\label{FigAttach}
\end{center}
\end{figure}

\pf We now define the reverse procedure. Given
$\sigma=\{B_1,B_2,\ldots,B_k\} \in \Pi_{[2,n]}$ and
$\omega=\omega_1\omega_2\ldots\omega_{k-1}\in [n]^{k-1}$, set
$\omega_k=1$. Define the inverse map $\phi_\sigma^{-1}: [n]^{k-1}
\rightarrow \mathcal{T}_{\succeq\sigma}$ as follows. For
$i=1,2,\ldots,k$:

\begin{enumerate}
\item Let $B(i)=\{b_1,b_2,\ldots,b_t\}_{<}$ be the largest block
of $\sigma \backslash \{B(1),B(2),\ldots,B(i-1)\}$ such that
$B(i)$ does not contain any number in $\{\omega_i,
\omega_{i+1},\ldots,\omega_{k-1}\}$.

\item Attach the vertices in $B(i)$ to $\omega_i$ according to the
following two cases:

{\it Case 1:} If $b_t>\omega_{i}$, we connect vertices
$b_1,b_2,\ldots,b_t$ and $\omega_i$ such that the edges between
$b_1,b_2,\ldots b_{t-1},\omega_i$ and $b_t$ are all pointing to
$b_t$ (See Figure \ref{FigAttach} (1));

{\it Case 2:} If $b_t<\omega_{i}$ we simply connect
$b_1,b_2,\ldots,b_t$ and $\omega_i$ such that all edges between
$b_1,b_2,\ldots,b_t$ and $\omega_i$ are all pointing to $\omega_i$
(See Figure \ref{FigAttach} (2)).

\end{enumerate}

It is easy to see that after all $k$ steps, we get a tree
$T:=\phi^{-1}_{\sigma}(\omega)\in \mathcal{T}_{\succeq \sigma}$.
One can easily check that $\phi_\sigma$ is a bijection.\qed

\begin{figure}[h,t]
\begin{center}
\setlength{\unitlength}{10pt}
\thicklines

\begin{picture}(35,15)
\put(0,0){
\begin{picture}(14,12)
\put(5,12){\circle*{0.3}}
\put(5.5,11.5){1}

\put(2,8){\circle*{0.3}} \put(2.5,7.5){7}
\put(5,12){\line(-3,-4){3}} \put(5,12){\vector(-3,-4){1.7}}
\put(2,4){\circle*{0.3}} \put(2.5,3.5){3}
\put(2,4){\line(0,1){4}}\put(2,4){\vector(0,1){2.4}}

\put(8,8){\circle*{0.3}} \put(8.5,7.5){9}
\put(5,12){\line(3,-4){3}} \put(5,12){\vector(3,-4){1.7}}
\put(5,4){\circle*{0.3}} \put(5.5,3.5){6}
\put(5,4){\line(3,4){3}} \put(5,4){\vector(3,4){1.7}}
\put(11,4){\circle*{0.3}} \put(11.5,3.5){5}
\put(11,4){\line(-3,4){3}}\put(11,4){\vector(-3,4){1.7}}

\put(8,0){\circle*{0.3}} \put(8.5,-0.5){2}
\put(8,0){\line(3,4){3}} \put(8,0){\vector(3,4){1.7}}
\put(11,0){\circle*{0.3}} \put(11.5,-0.5){4}
\put(11,0){\line(0,1){4}} \put(11,0){\vector(0,1){2.4}}

\put(14,0){\circle*{0.3}} \put(14.5,-0.5){8}
\put(11,4){\line(3,-4){3}} \put(11,4){\vector(3,-4){1.8}}
\end{picture}}


\put(20,0.5){
\begin{picture}(14,12)
\put(5,12){\circle*{0.3}} \put(5,12.3){\makebox(0,0)[b]{1}}

\put(8.5,8){\circle*{0.3}} \put(8.8,8){9}
\put(5.5,12){\line(3,-4){3}} \put(5.5,12){\vector(3,-4){1.7}} \put(7.5,10){$\bar{9}$}
\put(11.5,4){\circle*{0.3}} \put(11.8,4){5}
\put(11.5,4){\line(-3,4){3}} \put(11.5,4){\vector(-3,4){1.7}} \put(10.5,6){$\bar{5}$}

\put(5,3.7){\circle*{0.3}} \put(5.5,3.2){6}
\put(5,3.7){\line(3,4){3}} \put(5,3.7){\vector(3,4){1.7}} \put(5.5,5.7){$\bar{6}$}

\put(1.5,8){\circle*{0.3}} \put(2,7.5){7}
\put(4.5,12){\line(-3,-4){3}} \put(4.5,12){\vector(-3,-4){1.7}} \put(2,10){$\bar{7}$}

\put(1.5,3.5){\circle*{0.3}} \put(2,3){3}
\put(1.5,3.5){\line(0,1){4}}\put(1.5,3.5){\vector(0,1){2.4}}
\put(0.5,5){$\bar{3}$}

\put(8,-0.5){\circle*{0.3}} \put(8.5,-1){2}
\put(8,-0.5){\line(3,4){3}} \put(8,-0.5){\vector(3,4){1.7}} \put(8.5,1.5){$\bar{2}$}

\put(11,-0.5){\circle*{0.3}} \put(11.5,-1){4}
\put(11,-0.5){\line(0,1){4}}
\put(11,-0.5){\vector(0,1){2.4}}\put(11.3,0.5){$\bar{4}$}

\put(15,-0.5){\circle*{0.3}} \put(15.5,-1){8}
\put(12,3.5){\line(3,-4){3}} \put(12,3.5){\vector(3,-4){1.8}}\put(13.7,1.7){$\bar{8}$}

\end{picture}}
\end{picture}
\caption{A tree $T\in \mathcal{T}_{3221}$, with
$\phi(T)=8/569/37/24 \in \Pi_{3221}$, $\sigma=8/7/6/59/3/24
\prec\pi$, and $\phi_{\sigma}(T)=59715$.}\label{Figphisigma}
\end{center}
\end{figure}
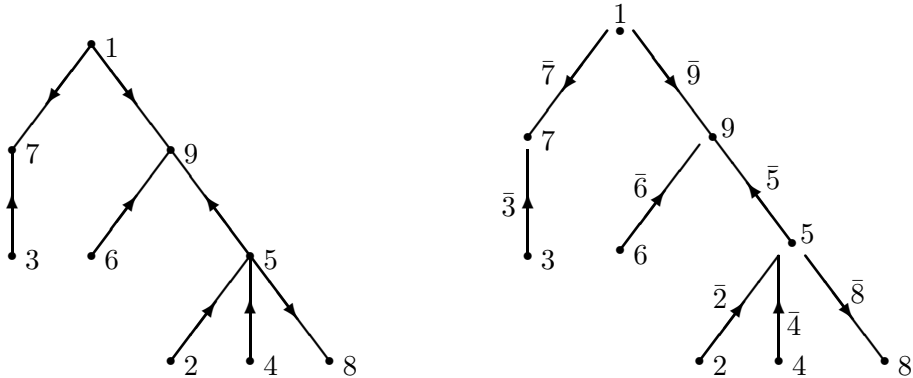

\begin{exam}
For the tree $T$ in Figure \ref{Figphisigma}, let
$\sigma=8/7/6/59/3/24$. We then have $B(1)=\{8\}$,
$\omega_{1}=c(B(1))=5$; $B(2)=\{6\}$, $\omega_{2}=c(B(2))=9$;
$B(3)=\{3\}$, $\omega_{3}=c(B(3))=7$; $B(4)=\{7\}$,
$\omega_{4}=c(B(4))=1$; $B(5)=\{2,4\}$, $\omega_{5}=c(B(5))=5$,
$B(6)=\{5,9\}$, $\omega_{6}=c(B(6))=1$ (which we don't write).
Thus we have $\phi_{8/7/6/59/3/24}(T)=59715\in [9]^{5}$.
\end{exam}

\noindent{\it Proof of Theorem \ref{main_thm}.} Let
$g(\sigma)=|\mathcal{T}_{\succeq\sigma}|$. From the bijection
$\phi_\sigma:\mathcal{T}_{\succeq\sigma}\rightarrow [n]^{
|\sigma|-1}$ we know that $g(\sigma)=n^{|\sigma|-1}$. Recall from
Section 2 that $f(\pi)=|\mathcal{T}_\pi|$. Since
$\mathcal{T}_{\succeq\sigma}=\bigcup_{\pi\succeq\sigma}\mathcal{T}_\pi$
is a disjoint union, we have
\begin{equation}\label{Mobius}
\sum_{\pi\succeq\sigma}f(\pi)=n^{k-1}, \mbox{ for any }\sigma
\in\Pi_{[2,n]}.
\end{equation}
It is now sufficient to prove that the unique solution of the
above equations is $f(\pi)=(n-1)!/(n-|\pi|)!$, for any $\pi \in
\Pi_{[2,n]}$.

First, since equation (\ref{Mobius}) holds for any
$\pi,\sigma\in\Pi_{[2,n]}$ such that $\pi\succeq\sigma$, we have,
by the poset structure of $\Pi_{[2,n]}$, that the solution $f$ to
equation (\ref{Mobius}) (valid for all $\sigma\in\Pi_{[2,n]}$) is
unique.

Second, let $\sigma=\{B_1, B_2, \ldots, B_k\}$. Then the interval
$[\sigma,\hat{1}_{[2,n]}]$ is isomorphic in an obvious way to the
lattice of partitions of the set $\{B_1, B_2, \ldots, B_k\}$.
Hence $[\pi,\hat{1}_{[2,n]}] \cong \Pi_{[k]}$, where
$\hat{1}_{[2,n]}$ is the maximum element of $\Pi_{[2,n]}$. Thus we
have
\begin{eqnarray*}
\sum_{\pi\succeq\sigma}\frac{(n-1)!}{(n-|\pi|)!}
&=&\sum_{\tau\in\Pi_{[k]}}\frac{(n-1)!}{(n-|\tau|)!}\\
&=&\sum_{j=1}^{k}S(k,j)\frac{(n-1)!}{(n-j)!}\\
&=&\frac{1}{n}\sum_{j=1}^{k}S(k,j)n(n-1)\cdots(n-j+1)\\
&=&n^{k-1},
\end{eqnarray*}
\noindent where $S(k,j)$ is the \emph{Stirling number of the
second kind}, i.e., the number of partitions of a $k$-set into $j$
blocks. The last equation follows from a standard Stirling number
identity, see e.g., identity (24d) in \cite[\S 1.4]{EC1}. Thus,
$(n-1)!/(n-|\pi|)!$ is a possible solution to equations
(\ref{Mobius}).

Hence, by uniqueness, we have $f(\pi)=(n-1)!/(n-|\pi|)!$.\qed

\noindent\emph{Remark: }In fact, given equation (\ref{Mobius}), we
can solve for $f$ by using the dual form of the M\"{o}bius
inversion formula:
\[f(\pi)=\sum_{\sigma \geq
\pi}\mu(\pi,\sigma)g(\sigma),\] where the coefficient
$\mu(\pi,\sigma)$ is the M\"{o}bius function of $\Pi_{[2,n]}$,
which can be calculated explicitly, \cite[Exam 3.10.4]{EC1}.

\section{The Real Bijective Map}

Although we gave a bijection $\phi_{\sigma}$ in Section 3, we
needed to prove Theorem \ref{main_thm} by solving equations, and
we still don't have a very good bijection that maps
$\mathcal{T}_\lambda$ to a set of cardinality $a_\lambda$ for any
$\lambda\vdash n-1$.

Let $\phi':\mathcal{T}_{\lambda}\rightarrow\Pi_\lambda\times
[n]^{k-1},T\mapsto (\pi,\phi_{\pi}(T))$, where $\pi=\phi(T)$.
Since $\phi_\pi$ is a bijection, we have that $\phi'$ is an
injection. Let $\Omega_\pi:=\phi_\pi(\mathcal{T}_\pi)$. Then
$\phi'(\mathcal{T}_{\lambda})=\{\{\pi\}\times
\Omega_\pi:\pi\in\Pi_\lambda\}=:(\Pi\times\Omega)_\lambda$. Thus,
$\phi':\mathcal{T}_{\lambda}\rightarrow (\Pi\times\Omega)_\lambda$
is the bijection we are looking for.

\begin{exam}
Assume $\pi=\{B_1,B_2\}$. For any $T\in\mathcal{T}_\pi$, we have
$\phi'(T)=(\pi,\max\{c(B_1),c(B_2)\})$, and $\Omega_\pi=[n-1]$,
$f(\pi)=n-1$.
\end{exam}

\begin{exam}
When $\lambda=\langle 1^{n-1}\rangle$, $\Pi_{\lambda}$ contains
only the partition $\hat{0}_{[2,n]}=n/n-1/\ldots/2$. As also
pointed out in the proof in Section 2,
$\mathcal{T}_{\hat{0}_{[2,n]}}$ is the set of all increasing trees
on $[n]$, in this case we have
$\Omega_{\hat{0}_{[2,n]}}=[n-1]\times[n-2]\times\cdots \times[1]$,
and for each $T\in \mathcal{T}_{\hat{0}_{[2,n]}}$ and
$\phi'(T)=(\pi,\omega)$, $\omega$ is the Pr\"{u}fer code of $T$.
\end{exam}

Though it seems quite hard to find what $\Omega_\pi$'s are, there
still exists a very good relation among them.

\begin{theo}\label{subsequence}
For any $\pi_1,\pi_2\in\Pi_{[2,n]}$, if $\pi_2\succ\pi_1$, we have
that for any $T\in\mathcal{T}_{\succeq\pi_2}$, $\phi_{\pi_2}(T)$
is a subsequence of $\phi_{\pi_1}(T)$. In particular, if
$\pi_2=\phi(T)$, we have $\phi'(T)=(\pi_2,\omega)$ and $\omega$ is
a subsequence of $\phi_{\pi_1}(T)$.
\end{theo}

\pf It suffices to prove the assertion for all covering pairs.
Assume that $\pi_2\cdot\!\!\!\succ\pi_1$. Thus there exist two
blocks $B$ and $B'$ of $\pi_1$ which become one block in $\pi_2$.

Assume
$\phi_{\pi_1}(T)=\omega=\omega_1\omega_2\cdots\omega_{k-1}$,
$\phi_{\pi_2}(T)=\omega^{\prime}=\omega'_1\omega'_2\cdots\omega'_{k-2}$
and $\omega_k=\omega'_{k-1}=1$. Then there exist $1\leq r,s\leq k$
such that $B$ and $B'$ are removed from $T$ at steps $r$ and $s$,
respectively, in process $\phi_{\pi_1}$. Assume, without loss of
generality, that $r<s$. Then it is easy to see that
$\omega'_l=\omega_l$ for $1\leq l< r$.

For step $r$ in process $\phi_{\pi_2}$, there are two cases:

\noindent \emph{Case 1:} If $B<B'$ and $\Star(B\cup B^{\prime})$
is a leaf-star, then it must be that $s=r+1$. At this step, we
remove $\Star(B\cup B^{\prime})$. Then $\omega'_r=\omega_s$,
$\omega'_l=\omega_l$ for $r<l\leq k-2$.

\noindent \emph{Case 2:} If $B>B'$, or $\Star(B\cup B^{\prime})$
is not a leaf-star, then we will remove $\Star(B\cup B^{\prime})$
at the step that we remove $\Star(B^{\prime})$ in the process
$\phi_{\pi_1}$, i.e., $\omega'_l=\omega_{l+1}$ for $r\leq l\leq
k-2$.

\noindent In both cases, we have that $\omega'$ is a subsequence
of $\omega$, i.e., $\phi_{\pi_2}(T)$ is a subsequence of
$\phi_{\pi_1}(T)$.\qed

\begin{exam}
Let $\pi=8/569/37/24$ and $\sigma=8/7/6/59/3/24$, so
$\pi\succ\sigma$. For the tree $T$ in Figure \ref{Figphisigma}, we
have $T \in \mathcal{T}_{\succeq \sigma}$, and
$\phi_{\pi}(T)=515$, which is a subsequence of
$\phi_{\sigma}(T)=59715$.
\end{exam}

By the proof of Theorem \ref{subsequence}, for any
$\sigma\in\Pi_{[2,n]}$ we can define a bijection from
$\bigcup_{\pi\succ\sigma}\Omega_\pi$ to
$[n]^{k-1}\backslash\Omega_\sigma$ such that each sequence will be
a subsequence of its image. Inductively using this bijection, we
can find out all $\Omega_\sigma$'s. But when $|\sigma|$ gets
larger and larger, it will become more and more difficult to find
out what this bijection is explicitly.

\section{Remarks}

We want to remark that the bijection we defined in Section 3 can
be considered as a generalization of the Pr\"{u}fer codes for
labelled trees: instead of deleting (attaching) vertices one by
one, we are dealing with groups of vertices with respect to a
partition of $[2,n]$. Moreover, the bijection $\phi^{\prime}$
together with Theorem \ref{subsequence} suggests a structure on
the set of labelled trees $\{\mathcal{T}_\pi:\pi\in\Pi_{[2,n]}\}$
as a lattice isomorphic to $\Pi_{[2,n]}$ under the map
$\mathcal{T}_\pi\mapsto\pi$.

The following problems are still interesting to consider.
\begin{enumerate}

\item{
Given $\pi\in\Pi_{[2,n]}$, Theorem \ref{subsequence} shows how to
find $\Omega_\pi$ explicitly, i.e., $\Omega_\pi$ is the subset of
$[n]^{|\pi|-1}$ with sequences corresponding to its subsequences
from $\Omega_{\sigma}$ deleted, for any $\sigma\succ\pi$. For
example, let $\pi=45/3/2$, we have:
\begin{equation*}
\Omega_\pi=\left\{\,
\begin{matrix}
  11 & 12 & \diagdown\!\!\!\!\!\!13\!\!\!\!\!\!\diagup & 14 & \diagdown\!\!\!\!\!\!15\!\!\!\!\!\!\diagup \\
  21 & 22 & \diagdown\!\!\!\!\!\!23\!\!\!\!\!\!\diagup & 24 & \diagdown\!\!\!\!\!\!25\!\!\!\!\!\!\diagup \\
  31 & 32 & \diagdown\!\!\!\!\!\!33\!\!\!\!\!\!\diagup & 34 & \diagdown\!\!\!\!\!\!35\!\!\!\!\!\!\diagup \\
  41 & 42 & 43 & \diagdown\!\!\!\!\!\!44\!\!\!\!\!\!\diagup & \diagdown\!\!\!\!\!\!45\!\!\!\!\!\!\diagup \\
  \diagdown\!\!\!\!\!\!51\!\!\!\!\!\!\diagup & \diagdown\!\!\!\!\!\!52\!\!\!\!\!\!\diagup & \diagdown\!\!\!\!\!\!53\!\!\!\!\!\!\diagup & \diagdown\!\!\!\!\!\!54\!\!\!\!\!\!\diagup & \diagdown\!\!\!\!\!\!55\!\!\!\!\!\!\diagup
\end{matrix}
\,\right\},
\end{equation*}
where $13,23,33,44$ correspond to its subsequences $1,2,3,4$ in
$\Omega_{45/23}$, $15,25,53,45$ correspond to its subsequences
$1,2,3,4$ in $\Omega_{3/245}$, $51,52,35,54$ correspond to its
subsequences $1,2,3,4$ in $\Omega_{345/2}$, and $55$ correspond to
its subsequence $\emptyset$ in $\Omega_{2345}$.

However, the ``corresponding relationship", between sequences and
its subsequences described inductively in the proof of Theorem
\ref{subsequence}, depends highly on the set
$\{\sigma\in\Pi_{[2,n]}:\sigma\succ\pi\}$, and it is not easy to
describe in general. Hence, it would be nice if one can give a
simple description of this relationship, and use it to characterize
$\Omega_\pi$.}

\item{In the proof of Theorem \ref{main_thm} in Section 2, we
mentioned that when $\lambda=\langle 1^{n-1}\rangle$, we can map
an increasing tree to a permutation of $[2,n]$ (\cite[\S
1.3]{EC1}). Is it possible to generalize this bijection to any
$\lambda$ by mapping a tree in $\mathcal{T}_\lambda$ to
$(\phi(T),w)$, where $w$ is a length $k-1$ permutation of
 an $(n-1)$-element set?}

\end{enumerate}

\vskip 2mm \noindent{\bf Acknowledgments.} The authors want to thank
Prof. Richard P. Stanley for his inspiring suggestions and helpful
discussions, and also thank the referee for helpful comments on an
earlier version of this paper. The first author is supported by the
National Science Foundation of China under Grant No.~10726048 and
No.~10801053, Shanghai Educational Development Foundation under the
Chenguang Project No.~2007CG29, and Shanghai Leading Academic
Discipline Project No.~B407.

\end{document}